\documentclass[11pt]{amsart}
\usepackage{amsfonts}

\newtheorem{theorem}{Theorem}[section]
\newtheorem{lemma}[theorem]{Lemma}
\newtheorem{proposition}[theorem]{Proposition}

\newtheorem{definition}[theorem]{Definition}

\newcommand{\C}{\mbox{$\mathbb{C}$}}
\newcommand{\N}{\mbox{$\mathbb{N}$}}

\newcommand{\K}{\mbox{$\mathbb{K}$}}

\newcommand{\A}{\mbox{${\mathcal A}$}}
\newcommand{\B}{\mbox{${\mathcal B}$}}
 
\begin{document}
\title[On Morita's Fundamental Theorem for 
C$^*-$algebras]{On Morita's Fundamental Theorem for C$^*-$algebras}   
 
\vspace{30 mm}
 
\author{David P. Blecher}
\date{May, 1997.  Revised June 1998.}
\thanks{* Supported by a grant from the NSF}
\thanks{The contents of this paper were announced at the joint 
meeting of the Canadian Operator Algebra Symposium, and the 
Great Plains Operator Theory Seminar, May 17-22, 1997}
\address{Department of Mathematics\\University of Houston\\Houston,
TX 77204-3476 }
\email{dblecher@@math.uh.edu}\maketitle
 
\maketitle
  
\vspace{10 mm}
 
\begin{abstract}
We give a solution, via operator spaces, of an old problem in the
Morita equivalence of C*-algebras.  Namely, we 
show that C*-algebras
are strongly Morita equivalent in the sense of Rieffel if and only if
their categories of left operator modules are isomorphic via completely
contractive functors.  Moreover, any such functor 
is completely isometrically
isomorphic to the Haagerup tensor product (= interior tensor product) with
a strong Morita equivalence bimodule.  
An operator module over a C$^*-$algebra
$\A$ is a closed subspace
of some B(H) which is left invariant under multiplication by $\pi(\A)$,
where $\pi$ is a *-representation of $\A$ on $H$.  The category 
$_{\A}HMOD$ of 
*-representations of $\A$ on Hilbert space
is a full subcategory of the category $_{\A}OMOD$ of operator modules.
Our main result remains
true with respect to subcategories of $OMOD$ which contain $HMOD$ and the
C$^*-$algebra itself.  It does not seem
possible to remove the operator space framework; 
in the very simplest cases there may exist no bounded equivalence functors
on categories with bounded module maps as
morphisms (as opposed to completely bounded ones).
Our proof
involves operator space techniques, together with
a C$^*-$algebra argument using compactness of
the quasistate space of a C$^*-$algebra, and lowersemicontinuity in
the enveloping von Neumann algebra.
\end{abstract}
 
\pagebreak
\newpage
 
\setcounter{section}{0}

\section{Notation, background and statement of the theorem}
In the early 70's M. Rieffel introduced and developed the notion of
strong Morita equivalence of C$^*-$algebras (see \cite{Ri2} for a 
good discussion and survey).  It has become a fundamental
tool in modern operator algebra and noncommutative geometry
(see \cite{Con} for example).  Briefly,
two C$^*-$algebras $\A$ and $\B$ are said to be
{\em strongly Morita equivalent} if
there is an $\A-\B-$bimodule $X$, which is a right C$^*-$module over
$\B$, and a left C$^*-$module over $\A$, such that the inner products
$_{\A}\langle \; \cdot  \; \vert \; \cdot \; \rangle$ and
$\langle \; \cdot  \; \vert \; \cdot \; \rangle_{\B}$ satisfy the relation
$_{\A}\langle \; x_1 \; \vert \; x_2 \; \rangle \; x_3
= x_1 \; \langle \; x_2 \; \vert \; x_3 \; \rangle_{\B}$, 
for $x_1,x_2,x_3 \in X$.  Also the span
of the range of these inner products must be norm dense in $\A$ and
$\B$ respectively.  Such $X$ is said to be
an $\A-\B-${\em strong Morita
equivalence bimodule}.

Our main result is a C$^*-$algebraic version of Morita's fundamental
theorem from pure algebra.
Namely,  we show that two C*-algebras
are strongly Morita equivalent if and only if
their categories of (left) operator modules are isomorphic via completely
contractive functors.  Moreover, 
any such functor is completely isometrically
isomorphic to the Haagerup tensor product (= interior tensor product) with
a strong Morita equivalence bimodule.  We use the context of operator 
spaces.  In previous papers \cite{BMP,Bna,Bhmo}
we showed that operator spaces, and more particularly operator 
modules, are an 
appropriate `metric' context for the C$^*-$algebraic theory of
strong Morita equivalence and the related theory of C$^*-$modules. 
Thus it was natural to look for a `fundamental Morita theorem' in this
category.        

Let us begin 
by establishing the common symbols and notations in this paper.  
We shall
use operator spaces quite extensively, and their connections to 
C$^*-$modules.  We refer the reader to \cite{Bna} and \cite{L2}
for missing background.
The algebraic background may be found in any account of Morita theory
for rings, such as \cite{AF}.  We have deliberately supressed some of the
purely algebraic calculations, since sentences consisting of
long strings of natural isomorphisms are not particularly interesting or
enlightening.  None of these supressed calculations are difficult,
and hopefully can be supplied without too much trouble by the reader.

We will use the symbols  
 ${\mathcal A}, {\mathcal B}$ for C$^*-$algebras; $a, b$ will
be generic elements of $\A$ and $\B$ respectively;
and $\{ e_\alpha \}, \{ f_\beta \}$ are contractive approximate identities
(c.a.i.'s) for $\A$ and $\B$ respectively.  We write $e(\A)$ for the 
enveloping von Neumann algebra of $\A$.   
 $H, K, \cdots$ are Hilbert spaces, $\zeta, \eta$ are
 typical elements in $H$ and
$K$ respectively, 
and $B(H)$ (resp. $B(H,K)$) is the space of bounded linear operators on 
$H$ (resp. from $H$ to $K$).  
We will reserve the symbols $Y$ and $Z$ for a right $\A$-module, or a left
$\B-$module, or an $\B-\A-$bimodule; it has generic element $y$ and $z$ 
respectively. 
Similarly, $X$ or $W$ will be a right $\B-$, left $\A-$, or
$\A-\B-$module, with generic element $x$ or $w$.  
     
Suppose that $\pi$ is a 
*-representation of $\A$ on 
Hilbert space $H$, and that $X$ is a closed subspace of 
$B(H)$ such that $\pi(\A) X \subset X$.  Then $X$ is a left $\A$-module.  
We shall assume that the module action is nondegenerate 
(= essential)\footnote{This means (for a left Banach module $X$ over 
$\A$, say)
that $\{\sum_{k=1}^n a_k x_k : n \in \N , \; a_k \in \A,
\; x_k \in X \}$ is
dense in $X$.  This is equivalent to saying that for any c.a.i.
 $\{ e_\alpha \}$ in $\A$,
$e_\alpha x \rightarrow x$ for all $x \in X$.}.
We say that such $X$, considered as
an abstract operator space and a left $\A$-module, is a
 left {\em operator module} over $\A$.  
By considering $X$ as an abstract operator 
space and module, we may forget about the particular $H, \pi$ 
used\footnote{It is sometimes useful, 
and equivalent, to allow $X$ in the definition
above, to be a subspace of $B(K,H)$, for a second Hilbert space $K$.}.
A theorem of Christensen-Effros-Sinclair \cite{CES} tells us that
the operator modules are exactly the operator spaces which are 
(nondegenerate) left $\A$-modules, such
that the module action 
is a `completely contractive' bilinear map (or equivalently,
the module action linearizes to a 
complete contraction $\A \otimes_h X \rightarrow X$,
where $\otimes_h$ is the Haagerup 
tensor product).  
We will use the facts
that submodules and  quotient modules
of operator modules, are again 
operator modules.  Also, if $X$ is a left operator module and
$E$ is an operator space, then the Haagerup tensor product
$X \otimes_h E$ is a left operator module.   This last fact follows
easily from the last definition of an operator module in terms of
the Haagerup tensor product, and the fact that that tensor product
is associative.   
We write $_{\A}OMOD$ for the category 
of left $\A$-operator modules.  The morphisms 
are $_{\A}CB(X,W)$, the {\em completely bounded} left 
$\A-$module maps.

We now turn to the category  $_{\A}HMOD$ 
of Hilbert spaces $H$ which are left $\A-$modules via a nondegenerate
$*-$representation of $\A$ on $H$
(denoted $_{\A}Hermod$ in \cite{Rieffel2}). 
If $H$ is a Hilbert 
space, and if $e_0$ is a fixed unit vector in $H$, then the space
of rank 1 operators $H^c = \{ \zeta \otimes e_0 \in B(H) : \zeta \in H \}$ is
clearly
an operator space, and indeed is clearly in $_{\A}OMOD$ if $\;
H \in \; _{\A}HMOD$.
As an operator space or operator module $H^c$  is independent
of the particular $e_0$ we picked.  
It is referred to in the literature as `Hilbert
column space'.
The $n$-dimensional Hilbert column space is written as $C_n$.
It is well known that for a linear map $T : H \rightarrow K$
between Hilbert spaces, the
usual norm equals the completely bounded norm of $T$
as a map $H^c \rightarrow K^c$.
Thus we see that the assignment 
$H \mapsto H^c$  embeds $_{\A}HMOD$ as a subcategory of 
$_{\A}OMOD$.  Henceforth we will view
it as a subcategory.
 
It is explained in \cite{Wi} that C$^*-$modules also 
possess canonical operator space structures, and so can be
viewed as objects in $OMOD$.  In  \cite{Bna} this idea is developed 
and, amongst other things, we showed that 
the well known interior tensor product of C$^*-$modules coincides
with their Haagerup tensor product as operator modules.
This fact is important in what follows.

If $X,W \in $ $_{\A}OMOD$ then 
$_{\A}CB(X,W)$ is an operator space
\cite{ERbimod}.     
In this paper we are concerned with functors between categories of
operator modules.  Such functors   
$F : $ $_{\A}OMOD \rightarrow $ $_{\B}OMOD$ are  
assumed to be  linear on spaces of morphisms.  Thus 
$T \mapsto F(T)$ from $_{\A}CB(X,W) \rightarrow
$ $_{\B}CB(F(X),F(W))$ is linear, 
for all pairs of objects $X,W \in $ $_{\A}OMOD$.  
We say $F$ is  {\em completely contractive}, if
this map $T \mapsto F(T)$
is completely contractive, 
for all pairs of objects $X,W \in $ $_{\A}OMOD$.
We say two functors  $F_1, F_2 : $ $_{\A}OMOD \rightarrow $ $_{\B}OMOD$ 
are
(naturally) completely isometrically isomorphic, if they are 
naturally isomorphic in
the sense of category theory 
\cite{AF}, with the natural transformations being 
complete isometries.  In this case we write 
$F_1 \cong F_2$ {\em completely isometrically}.  
 
\begin{definition}
\label{defome} We say that two C$^*-$algebras 
$\A$ and $\B$ are {\em operator Morita equivalent}
if there exist completely contractive  
functors $F : $ $_{\A}OMOD \rightarrow $ $_{\B}OMOD$
and $G : $ $_{\B}OMOD \rightarrow $ $_{\A}OMOD$, such that $F G \cong Id$ and 
$G F \cong Id$ completely isometrically.  
Such  $F$ and $G$ will be called {\em operator 
equivalence functors}.
\end{definition}    

We can now state our main theorem.  Its proof, which occupies  \S 2 and 3,
involves operator space techniques, together with
a C$^*-$algebra argument using compactness of
the quasistate space $Q$ of a C$^*-$algebra, and lowersemicontinuity in
the enveloping von Neumann algebra.

\begin{theorem}
\label{FT}  
Two C$^*-$algebras $\A$ and $\B$
are strongly Morita equivalent if and only if they
are operator Morita equivalent.  Suppose that
 $F, G$ are the operator equivalence 
functors, and set $Y = F(\A)$ and $X = G(\B)$.  Then $X$ is an
$\A-\B-$strong Morita equivalence bimodule,  $Y$ is a
$\B-\A-$strong Morita equivalence bimodules, and  $Y$  is 
unitarily equivalent 
to the conjugate C$^*-$bimodule $\bar{X}$ of $X$.  Moreover, 
$F(W) \cong Y \otimes_{h\A} W \cong $ $_{\A}\K(X,W)$  
completely isometrically 
isomorphically (as $\B-$operator modules), for all $W \in $ $_{\A}OMOD$.  
Thus $F \cong Y \otimes_{h\A} - \; \cong $ 
$_{\A}\K(X,-)$ completely isometrically.   
Similarly $G \cong X \otimes_{h\B} - \; \cong $ $_{\B}\K(Y,-)$ 
completely isometrically.   Also $F$ maps the subcategory
$_{\A}HMOD$ to $_{\B}HMOD$, and
the subcategory of C$^*$-modules to the C$^*-$modules (on which 
subcategories the Haagerup tensor product above coincides with
the interior tensor product).  Similar statements hold for $G$.  
\end{theorem}  

\vspace{6 mm}

We remind the reader that $_{\A}\K(X,W)$ was defined 
in \cite{Bna} to be the
norm closure in $_{\A}CB(X,W)$ of the span of the rank one operators
$\langle \cdot \; \vert \; x \; \rangle w$, for $x \in X, w \in W$.  
The symbol $\otimes_{h\A}$ denotes the module Haagerup tensor
product over $\A$ (see \cite{BMP} or \cite{Bna}).  

\vspace{6 mm}

\noindent {\bf Remark 1).}  The one direction of the ``if and only if''
of the theorem is easy and was    
noted in \cite{BMP}.  For completeness we sketch the short
argument here.  Namely, if
$X$ is a strong 
Morita equivalence bimodule for a strong Morita equivalence of
$\A$
and $\B$, and if $Y = \bar{X}$ is the conjugate C$^*-$module, then define
$F(W) = Y \otimes_{h\A} W, $ and
$G(Z) = X \otimes_{h\B} Z$.  Since the Haagerup
tensor product is functorial, $F$ and $G$ are functors.
By the associativity of
the module Haagerup tensor product, and the fact that
this tensor product equals the
interior tensor product where the
latter is defined, we obtain that
$$
GF(W) \cong X \otimes_{h\B} (
Y \otimes_{h\A} W) \cong 
 (X \otimes_{h\B} Y) \otimes_{h\A} W
\cong \; \A \otimes_{h\A} W \cong W
$$
completely isometrically, and as $\A-$modules.  Similarly $FG \cong Id$
completely isometrically.  So $\A$ and $\B$ are operator Morita equivalent.

\vspace{6 mm}

\noindent {\bf Remark 2).}   One can adapt the statement of the 
theorem to allow the operator
equivalence functors to be defined on not all of $OMOD$,
but only on a subcategory ${\bf D}$ of
$OMOD$ which contains $HMOD$ and the C$^*-$algebra itself.
Our proof goes through verbatim, 
except that for the part in \S 2 that equivalence functors 
preserve $HMOD$.  For this part to work, the subcategory 
${\bf D}$ should be closed under two or three
operations which we leave to the interested
reader  to abstract.

We also remark that the proof would become a little simpler
if we are willing to assume that the functors
concerned are `strongly continuous' (by which we mean that
$F(T_\lambda)$ converges
point norm to $F(T)$ whenever $T_\lambda$ is a bounded net in
$_{\A}CB(X,W)$ converging point norm to $T$).
This argument, which was in the original version of this paper,
has been omitted for the sake of brevity.

\vspace{6 mm}

\noindent {\bf Remark 3).}  
The reader may question the necessity of using operator
spaces, and completely contractive or completely isometric maps and functors. 
However it is not too hard to show that 
even in the very simplest case, where $\A = \C$,
$\B = M_n$ (the $n \times n$ scalar 
matrices), and if we write ${\bf D}$ for either
the category 
of left Banach modules, or the category of operator modules but
with bounded module 
maps as opposed to completely bounded ones, that there exists no
isometric equivalence functor 
 $F : $ $_{\A}{\bf D} \rightarrow $ $_{\B}{\bf D}$.         
In these categories there are too many morphisms; one needs to restrict 
attention to the completely bounded ones.  If one replaces $\B$ by
the compact operators on $\ell^2$, 
there exists no bounded equivalence functor
(see also \cite{Gk}).

Indeed one runs into problems using bounded module maps as morphisms
if one picks the smallest categories containing $HMOD$ and the algebra
itself.  Namely, suppose that $\A$ and $\B$ are strongly 
Morita equivalent, with $\A-\B-$equivalence bimodule $X$ and
dual bimodule $Y \cong \bar{X}$.  Let $_{\A}{\bf C}$ be the category of 
left Banach (or operator)
$\A-$modules consisting of $_{\A}HMOD$, $\A$ and $X$ (the
latter two viewed as left $\A-$modules).   Let $_{\B}{\bf D}$
consist of $_{\B}HMOD$, $B$ and $Y$.  Morphisms in both 
categories are the bounded module maps.  Take $F$ to be the obvious
functor,
namely the one that maps $\A$ to $Y$, $X$ to $B$, and on $_{\A}HMOD$
is the interior tensor product with $Y$.   Define $G : $ $_{\B}{\bf D}
\rightarrow $ $_{\A}{\bf C}$ similarly.  Again it is easy to 
check that even in the simplest cases $F$ and $G$ are not necessarily
contractive or bounded.

\vspace{4 mm}

\noindent {\bf Remark 4).}
 W. Beer proved in \cite{Beer}
that two unital C$^*-$algebras are 
strongly Morita equivalent if and only if they are algebraically
Morita equivalent.  Our theorem may be viewed as an extension to the 
general case which also has the advantage of characterizing the
equivalence functors up to (complete) isometry.  Also,
in Beer's theorem
one produces the C$^*-$module by finding a similarity of
 an idempotent in a matrix algebra to a selfadjoint idempotent,
whereas our equivalence bimodule comes directly from the functor.       

In \cite{Bna} we gave another C$^*-$algebraic analogue of
Morita's fundamental theorem in terms
of categories of C$^*-$modules; but
that theorem was much less satisfying.
The definition of a C$^*-$module is not too far from
that of a strong Morita equivalence, so that while that theorem
was not quite tautological, it was certainly not very
deep\footnote{Indeed the proof
of the aforementioned
theorem in \cite{Bna} is rather too long: as we noted
in the galley proofs to
that paper, G. Skandalis has shown us a shorter proof.}.
It seems much more surprizing, at least to us, that strong Morita
equivalence should be related to the category of operator modules.
After all, the notion of an operator module has nothing to do
with the notion of strong Morita equivalence.
Another `drawback' of the theorem in \cite{Bna} is that
the category of C$^*-$modules does not contain the category
$_{\A}HMOD$ of Hilbert space modules.

\vspace{5 mm}

\setcounter{section}{1}

\section{Preliminary Lemmas}

Throughout this section $\A$ and $\B$ are C$^*-$algebras, and
$F : $ $_{\A}OMOD \rightarrow $ $_{\B}OMOD$ 
is an operator equivalence functor, with `inverse' $G$ (see Definition
\ref{defome}).  We set $Y = F(\A)$ and $X = G(\B)$.  
For a a left module $W$  over $\A$, say, and $w \in W$,
we write $r_w$ for the map from $\A \rightarrow
W$ which is simply right multiplication by $w$.

\begin{lemma}
\label{coron}
Let  $W \in $ $ _{\A}OMOD$.
Then $w \mapsto r_w$ is a complete
isometry of $W$ into $_{\A}CB(\A,W)$.  Indeed,
 $W$ is completely 
isometrically isomorphic to $\{ T \in $ $_{\A}CB(\A,W) : T \; r_{e_\alpha}
\rightarrow T$ in norm $\; \}$, where $\{ e_\alpha \}$ is a c.a.i. for 
$\A$.  If $W$ is also 
a Hilbert space, then the map above is a completely isometric
isomorphism  $W \cong $ $_{\A}CB(\A,W)$.
\end{lemma}

This is a simple consequence of the existence of a
c.a.i. in any C$^*-$algebra.
The following lemma will be used extensively without comment.  It's
proof is just as in pure algebra (\cite{AF} Proposition 21.2).

\begin{lemma}
\label{bf}  If $V, W \in $ $_{\A}OMOD$ then the map $T \mapsto F(T)$ gives
a completely isometric
surjective linear isomorphism $_{\A}CB(V,W) \cong $ $_{\B}CB(F(V),F(W))$.
If $V = W$ this map is a completely isometric
isomorphism of algebras.
\end{lemma}

If $E$ is an operator space, then the space $M_{m,n}(E)$ of $m \times n$
matrices with entries in $E$, is also an operator space in a canonical 
way.  We write $C_m(E)$ and $R_m(E)$ for the operator spaces
$M_{m,1}(E)$ and $M_{1,m}(E)$.  If $W \in $ $_{\A}OMOD$, then it is easy to 
see that $R_m(W)$ and $C_m(W)$ are again in  $_{\A}OMOD$.   

For $n
 = 1, \cdots , m$,  write
$i_n$ (resp. $\pi_n$) for the canonical coordinatewise
inclusion (resp. projection) map of $W$ into
the direct sum $C_m(W)$ or $R_m(W)$ (resp. from the direct sum onto $W$).   
Then $\pi_n i_k = \delta_{k,n} Id_W$ for each $n,k$ (where
$\delta_{k,n}$ is the Kronecker delta), and $\sum_n i_n \pi_n = Id$.  
Applying the functor $F$ gives maps $F(i_n) : F(W) \rightarrow F(R_m(W))$,
and $F(\pi_n) : F(R_m(W)) \rightarrow F(W)$, with $F(\pi_n) F(i_k) = 
\delta_{k,n} Id$
for each $n, k$, and $\sum_n F(i_n) F(\pi_n) = Id$.  These formulae yield 
a canonical algebraic isomorphism $F(R_m(W)) \cong R_m(F(W))$.  Similarly
in the $C_m(W)$ case.  We now prove these isomorphisms are 
completely isometric: 

\begin{lemma}
\label{dis}  For any $W \in $ $_{\A}OMOD$, we have
$F(R_m(W)) \cong R_m(F(W))$ and $F(C_m(W)) \cong C_m(F(W))$
completely isometrically isomorphically.
\end{lemma}

\begin{proof}  
In the $R_m(W)$ case, note $[\pi_1 , \cdots , \pi_m] \in
 R_m(_{\A}CB(R_m(W),W))$, and
it has norm 1 (as may be seen by noting that it corresponds to the
identity map after
employing the canonical completely isometric identification
$R_m(_{\A}CB(R_m(W),W)) \cong $ $_{\A}CB(R_m(W),R_m(W))$).
Applying $F$, we find $J = [F(\pi_1) , \cdots , F(\pi_m)]
\in R_m(_{\B}CB(F(R_m(W)),F(W)))$ has norm 1.  However, via 
the canonical completely isometric  
isomorphism of $R_m(_{\B}CB(F(R_m(W)),F(W)))$ with 
$_{\B}CB(F(R_m(W)),R_m(F(W)))$, $J$ corresponds to the canonical morphism 
$F(R_m(W)) \rightarrow R_m(F(W))$.  So this latter morphism is a
complete contraction.  Similarly the 
canonical morphism $G(R_m(F(W))) \rightarrow R_m(GF(W)) \cong R_m(W)$
is a complete contraction.  Applying $F$ to this morphism, gives a
complete contraction $FG(R_m(F(W))) \rightarrow F(R_m(W))$, which yields a
complete contraction $R_m(F(W)) \rightarrow F(R_m(W))$.   This 
proves the lemma for $R_m(W)$.  The $C_m(W)$ case is similar.  
\end{proof}

\vspace{8 mm}
 
In the remainder of
this section we show that $F$ takes the subcategory $_{\A}HMOD$
to $_{\B}HMOD$, and similarly for $G$.   
Choose $H \in $ $_{\A}HMOD$, and recall
that  $H$ may be identified with $H^c \in $ $_{\A}OMOD$.   We will 
 show that $F(H^c) \in $ $_{\B}HMOD$, or equivalently, that
$F(H^c)$ is a column Hilbert space.  For this we need
the following functorial characterization of column Hilbert 
space:  

\begin{proposition}
\label{Pis} 
Let $E$ be an operator space.  Then $E$ is completely isometrically
isomorphic to a Hilbert column space if and only if the identity 
map $E \otimes_{min} C_m \rightarrow E \otimes_h C_m$ is a 
complete contraction for all $m \in \N$.  
\end{proposition}

\begin{proof}
The ($\implies$) direction is easy and is omitted \cite{B2,ERsd}.
A simple proof of the other direction
may be found in 
\cite{ERgp} (Theorem 4.1, setting $q=1$).   
For completeness, we 
sketch a slight simplification of their argument.  We use
canonical operator space identifications, which may be found in
\cite{BP1,ERsd,B2}, and the notation of \cite{B2}.  
By the complete injectivity of the minimal and Haagerup tensor product,
 (see \cite{BP1} for example),
and the fact that column Hilbert space is determined by its finite 
dimensional subspaces being column space,
it follows that $E \otimes_{min} H^c =
E \otimes_h H^c$, for any Hilbert space $H$.  Choose $H$ 
so that $E \subset B(H)$.
The last ``$=$" may be 
rewritten as $H^c \otimes_h E = H^c \overset{\frown}{\otimes} E$,
where $\overset{\frown}{\otimes}$ is the 
operator space projective tensor product.   
Next, recall that  the functors $H^r \otimes_h -$ and 
$H^r \overset{\frown}{\otimes} -$ are the same. 
Applying this functor to the identity
$H^c \otimes_h E = H^c \overset{\frown}{\otimes} E$
yields the identity
$S_1(H) \otimes_h E = S_1(H) \overset{\frown}{\otimes}
 E$, where $S_1(H)$ is the 
operator space predual of $B(H)$.   Taking the operator 
space dual yields $CB(E,B(H)) = \Gamma^c(E,B(H))$.  Thus the
inclusion map  $E \subset B(H)$ factors
through Hilbert column space.  Hence $E$ is Hilbert column space.
\end{proof}  

\vspace{7mm}  

We remark that Pisier has shown us that
the last result is true with the word ``complete''
removed.

To use this to prove that $K = F(H^c)$ is a column Hilbert space,
we first remind the reader that for Hilbert column
spaces, all operator space tensor norms coincide \cite{ERsd,B2},
thus $C_m(H^c) \cong H^c \otimes_{min} C_m \cong 
H^c \otimes_h C_m$ completely isometrically.
Using this and Lemma \ref{dis}  we have (completely 
isometrically):  
$$K \otimes_{min} C_m \cong C_m(F(H^c)) \cong F(C_m(H^c)) \cong
F(H^c \otimes_h C_m) \cong F(G(K) \otimes_h C_m) \eqno{(*)}   
$$
since $G(K) \cong H^c$.  Next we remark that
 there is a canonical complete contraction
$G(K) \otimes_h C_m \rightarrow G(K \otimes_h C_m)$.  To explain this
map, first consider the map $G(K) \rightarrow  $ $_{\B}CB(Y,K)$
given by the
following sequence of maps: 
$$
G(K) \; \rightarrow \; \; _{\A}CB(\A,G(K))
\; \cong \; \; _{\B}CB(Y,FG(K)) \; \cong \; \; _{\B}CB(Y,K) . 
\eqno{(**)}
$$

The $\rightarrow$ in (**) comes from Lemma \ref{coron},
whereas the $\cong$ comes from applying the equivalence functor
(see Lemma \ref{bf}). 
Using (**) we get a sequence of completely contractive
module maps:
$$ 
G(K) \otimes_h C_m \rightarrow \; _{\B}CB(Y,K) \otimes_h C_m
\rightarrow \; _{\B}CB(Y,K \otimes_h C_m) \cong  
\; _{\A}CB(\A,G(K \otimes_h C_m))  . \eqno{(***)} 
$$

The second $\rightarrow$ in (***) comes about because
any $T \in CB(Y,K)$ and $z \in C_m$ gives a map 
in $CB(Y,K \otimes_h C_m)$ given by $y \mapsto T(y) \otimes z$.  
Moreover it is easy to check that
 this prescription gives a complete contraction 
$_{\B}CB(Y,K) \otimes_h C_m
\rightarrow $ $_{\B}CB(Y,K \otimes_h C_m)$.
The $\cong$ in (***) comes from applying the equivalence 
functor.

If one checks through (***) one finds that the composition 
of the maps lands up in
$G(K \otimes_h C_m)$ inside $_{\A}CB(\A,G(K \otimes_h C_m))$.  That is,
(***) gives a map
 $G(K) \otimes_h C_m \rightarrow G(K \otimes_h C_m)$.
Applying $F$ to this last map and putting this together with (*)
  gives us a 
complete contraction

$$
K \otimes_{min} C_m \cong F(G(K) \otimes_h C_m)
\rightarrow F(G(K \otimes_h C_m)) \cong K \otimes_h C_m
$$
Thus we have obtained a complete contraction $K \otimes_{min} C_m
\rightarrow K \otimes_h C_m$ which, one can easily check,
up to complete isometry, is
the canonical map between these spaces.  Appealing to 
Proposition \ref{Pis}
completes the argument of this section.   

\setcounter{section}{2}

\section{Completion of the proof of the main theorem} 

Again $\A, \B, F, G, X, Y$ are as in the previous section, but now
we fix $H \in $ $_{\A}HMOD$ to be the Hilbert space of the universal 
representation 
of $\A$, and fix $K = F(H)$.  Then $e(\A) \subset B(H)$, where
$e(\A)$ is the enveloping von 
Neumann algebra of $\A$.  By \S 2, $F$ and $G$
restrict to an equivalence of $_{\A}HMOD$ with $_{\B}HMOD$. 
By elementary C$^*-$algebra facts,
$F$ and $G$ restricted to $HMOD$ are automatically normal
*-functors in the sense of \cite{Rieffel2}.
By 
\cite{Rieffel2} Propositions 1.1, 1.3 and 1.6, $\B$ acts faithfully
on $K$, and if we regard $\B$ as a subset of $B(K)$ then
the weak operator closure
$\B''$ of $\B$ in $B(K)$, is W*-isomorphic to $e(\B)$.  We shall indeed 
regard $\B$ henceforth as a subalgebra of $B(K)$.  
We shall 
need the fact, from \cite{Rieffel2} Proposition 4.9,
that if $H^\infty$ is the 
Hilbert space direct sum of a countably infinite number of copies
of $H$, then $F(H^\infty) \cong K^\infty$, and similarly
$G(K^\infty) \cong H^\infty$ .

It is important in what follows to keep in mind the canonical 
right module action of $\B$ on $X$.
$x b = F(r_b)(x)$, for $x \in X , b \in \B$, where
as in the previous section $r_b : \B \rightarrow \B 
: c \mapsto cb$.
Similarly, $Y$ is canonically a $\B-\A-$bimodule.

There is a left $\B-$module map
$Y \otimes X \rightarrow F(X)$ defined by
$y \otimes x \mapsto F(r_x)(y)$.  Since $F(X) = FG(\B) \cong \B$, 
we get a 
left $\B$-module map
$Y \otimes X \rightarrow \B$, which we shall write as
$[\cdot]$.  Simple algebra shows that $[\cdot]$ is a
$\B-\B-$bimodule map, but this will not be explicitly needed.   
 In a similar way we get a module map 
$(\cdot) : X \otimes Y \rightarrow \A$.  
In what
follows we may use the same notations for the `unlinearized' bilinear
maps, so for example we may use the symbols $[y,x]$ for $[y \otimes x]$.
We now show that these maps have dense range.   
By way of contradiction, 
suppose that the closure of the range of $[\cdot]$ 
is a proper submodule $I$ of $\B$.
Let $Z = \B / I$, regarded as a left $\B-$operator module (see Lemma 2.1
in \cite{BMP}), and let $\pi :$ $\B \rightarrow Z$ be the nonzero 
quotient map.  Then $G(\pi) : X \rightarrow G(Z)$ is nonzero.  So there
exists $x \in X$ such that  $G(\pi) r_x \neq 0$.
Applying
$F$ we obtain $FG(\pi) F(r_x) \neq 0$, so that for some $y \in Y$,
$FG(\pi) F(r_x) (y) \neq 0$.  By the definition of $[\cdot]$ this implies
that $\pi([y \otimes x]) \neq 0$, 
which contradicts the definition of $\pi$.
Thus $[\cdot]$ (and similarly $(\cdot)$) has dense range. 

It should be pointed out that if we are attempting to prove the main
theorem, but with $OMOD$ replaced by a subcategory (as discussed in
Remark 2 in  \S 1), then the 
argument of the last paragraph seems to require 
that the subcategory be closed under certain quotients.  However, the
last paragraph can be replaced by an argument  which avoids a quotient
in the subcategory.  Namely, pick a faithful (nonzero) representation
of $\B / I$ on a Hilbert space $K$ say.  Then $K$ can be regarded in
a canonical way as an object in $_{\B}HMOD$.  Then there is a
nonzero morphism $S$ from $\B / I$ to $K$.  Replace the map 
$\pi$ in the previous paragraph by $S \circ \pi$, 
which is a nonzero morphism from $\B$ to $K$, and proceed in the same way.

\begin{lemma}
\label{xin}  The canonical maps $X \rightarrow $ $_{\B}CB(Y,B)$ and
$Y \rightarrow $ $_{\A}CB(X,A)$ induced by $[\cdot]$ and $(\cdot)$
respectively, are complete isometries.
\end{lemma}

\begin{proof}
Using Lemmas \ref{coron}, \ref{bf}, and the fact that $F(X) =
FG(\B) \cong \B$, we have
$X \subset $ 
$_{\A}CB(\A,X) \cong $ $_{\B}CB(Y,F(X)) \cong $ $_{\B}CB(Y,\B)$
completely isometrically.  Sorting through these identifications shows
that an element  $x \in X$ corresponds to the map $y \mapsto [y,x]$ in
$_{\B}CB(Y,\B)$.  A similar proof works for $(\cdot)$.
\end{proof}

\vspace{7 mm}

The following maps $\Phi : Y  \rightarrow B(H,K)$, and $\Psi : X
\rightarrow B(K,H)$ will play a central role in the remainder of the 
proof.  Namely, $\Phi(y)(\zeta) = F(r_\zeta)(y)$, and 
$\Psi(x)(\eta) = \omega_H 
G(r_\eta)(x)$, where $\omega_H : GF(H) \rightarrow
H$ is the $\A-$module map
coming from the natural transformation $GF \cong Id$.  Since
$\omega_H$ is an isometric surjection between Hilbert space it is unitary,
which will be important below.  It
is straightforward algebra to check that:
$$
\Psi(x) \Phi(y) = (x,y) \hspace{15mm} \& \hspace{15mm}
 \Phi(y) \Psi(x) = [y,x] V \; \;
  \eqno{(1)} 
$$
for all $x \in X, y \in Y$, and $V \in B(K)$ is a unitary operator in 
$\B'$ composed of two natural transformations.  
The $V$ will not play a significant role, since we will mostly be working
with expressions such as $[y,x]^*[y,x]$ which by the above, and since 
$V$ is unitary and in $\B'$, equals $\Psi(x)^*\Phi(y)^*\Phi(y)\Psi(x)$. 

Before we begin the next lemma, we remark that for any Hilbert
spaces $H,K$, since $CB(H^c,K^c) = B(H,K)$ 
completely isometrically (see \cite{Wi,ERsd,B2}),  the norm of a matrix
$[T_{ij}] \in M_n(B(H,K))$ can be calculated by the formula: 
$$ 
\Vert [T_{ij}] \Vert = \sup \{ \Vert [T_{ij}(\zeta_{kl})] \Vert
: [\zeta_{kl}] \in Ball(M_m(H^c)) , m \in \N \}  \eqno{(2)}
$$

\begin{lemma}
\label{phci}
The map $\Phi$ (resp. $\Psi$) is a completely isometric $\B-\A-$module
map (resp. $\A-\B-$module map).   Moreover, $\Phi(y_1)^*\Phi(y_2) \in
\A'' = e(\A)$ for all $y_1, y_2 \in Y$, and $\Psi(x_1)^*\Psi(x_2) \in
\B''$ for $x_1, x_2 \in X$.
\end{lemma}
 
\begin{proof}
We shall simply prove the assertions for $\Phi$; those for $\Psi$
are similar.  The module map assertions are fairly clear, for instance
$\Phi(ya)(\zeta) = F(r_\zeta)(ya) = F(r_\zeta)F(r_a)(y) = 
F(r_{a\zeta})(y) = \Phi(y)(a\zeta)$ .  Next we show the $\A''$ assertion.
By Lemma \ref{bf}, we have a C$^*-$isomorphism $T \mapsto F(T) : $ $\A'
 = $ $_{\A}CB(H^c) \rightarrow \B' = $  $_{\B}CB(K^c)$. Note   
$\Phi(y)T(\zeta) = F(r_{T(\zeta)})(y) = F(T) F(r_\zeta)(y) = F(T) \Phi(y)
(\zeta), $ for $T \in $ $\A'$, and so also $(T \Phi(y)^*)^*
= \Phi(y) T^* = F(T^*) \Phi(y) = (\Phi(y)^* F(T))^*$.  
Together these imply that $\Phi(y_1)^*\Phi(y_2) \in $ $\A''$.  
The matching assertion for $\Psi$ has the additional complication
of the maps $\omega_H$, however since they are unitary as remarked above, 
they disappear from the calculation. 
Finally, we turn to the complete isometry.  The equalities in the following 
calculation follow from, in turn, formula (2) above, 
the definition of $\Phi$, Lemma \ref{coron}, Lemma \ref{bf}, the definition
of $(\cdot)$, formula (2) again, and Lemma \ref{xin}:
$$
\begin{array}{ccl}
\Vert [\Phi(y_{ij})] \Vert & = & \sup 
\{ \Vert [\Phi(y_{ij})(\zeta_{kl})] \Vert
: [\zeta_{kl}] \in Ball(M_m(H^c)) , m \in \N \} \\
& = & \sup \{ \Vert [F(r_{\zeta_{kl}})(y_{ij})] \Vert
: [\zeta_{kl}] \in Ball(M_m(H^c)) , m \in \N \} \\
& = & \sup \{ \Vert [F(r_{\zeta_{kl}}) \; r_{y_{ij}} ] \Vert :
 [\zeta_{kl}] \in Ball(M_m(H^c)) , m \in \N \} \\
& = & \sup 
\{ \Vert [GF(r_{\zeta_{kl}}) \; G(r_{y_{ij}}) ] \Vert: [\zeta_{kl}] 
\in Ball(M_m(H^c)) , m \in \N \} \\
& = & \sup \{ \Vert [GF(r_{\zeta_{kl}}) \; G(r_{y_{ij}}) (x_{pq}) ] \Vert:
[\zeta_{kl}] \in Ball(M_m(H^c)) , [x_{pq}] \in Ball(M_r(X)) \} \\ 
& = & \sup \{ \Vert [(x_{pq},y_{ij})\zeta_{kl}] \Vert :
[\zeta_{kl}] \in Ball(M_m(H^c)) , [x_{pq}] \in Ball(M_r(X)) 
, m,r \in \N \} \\   
& = & \sup \{ \Vert [(x_{pq},y_{ij})] \Vert : [x_{pq}] \in Ball(M_r(X)) ,
r \in \N \} \\
& = & \Vert [y_{ij}] \Vert
\end{array} 
$$
Thus $\Phi$ is a complete isometry.
\end{proof}  

We now proceed towards showing that
\begin{theorem}
\label{slide}  Suppose that 
$\Psi(x)^*\Psi(x)$, which is in $\B''$ by the previous lemma,
is actually in $\B$
for all $x \in X$; and suppose that 
$\Phi(y)^*\Phi(y) \in $ $\A$ for all $y \in Y$.  
Then all the conclusions of our main theorem hold.
\end{theorem}

\begin{proof}
If $\Psi(x)^*\Psi(x) \in \B$ for all $x \in X$, then by
the polarization
identity, and the previous
lemma,  $X$ is a RIGHT C$^*-$module over $\B$ with i.p. 
$\langle \; x_1 \; \vert \; x_2 \; \rangle_{\B} = \Psi(x_1)^*\Psi(x_2)$ .  
We can also 
deduce that  $X$ is a LEFT C$^*-$module over $\A$
by setting 
$_{\A}\langle \; x_1 \; \vert \; x_2 \; \rangle = \Psi(x_1) \Psi(x_2)^*$.
This last quantity may be seen to lie in $\A$ by using the
polarization identity and the following argument:
Since the range of $(\cdot)$ is dense in $\A$, we can find a
c.a.i. $\{ e_\alpha \}$ for $\A$, with terms of the form $e_\alpha
= \sum_{k=1}^{n}
(x_k,y_k) = \sum_{k=1}^n \Psi(x_k) \Phi(y_k)$ (using equation
(1)).  Here $n, x_k, y_k$ depend on $\alpha$.  Then $\{ e^*_\alpha \}$ is
also a c.a.i. for $\A$.    Since $\Psi(x)^* =
\lim_\alpha \; \Psi(e_{\alpha}^* x)^* = \lim_\alpha \;
 \Psi(x)^* e_\alpha$, it follows that
$\Psi(x) \Psi(x)^*$  is a norm limit of finite sums of terms of the
form $\Psi(x) (\Psi(x)^* \Psi(x_k)) \Phi(y_k) 
= \Psi(x) \Phi(b y_k)  = (x,by_k) \in $ $\A$, where
 $b = \Psi(x)^*\Psi(x_k)
\in $ $\B$.  Thus $\Psi(x) \Psi(x)^* \in $ $\A$.  

A similar argument shows that $Y$ (or equivalently $\Phi(Y)$) is 
both a left and right C$^*-$module.  At this point we can therefore
say that
the right module actions on $Y$ and $X$ are nondegenerate.
Notice also, that
if we choose a contractive approximate identity for $\A$ of
form $e_\alpha = \sum_k \Psi(x_k) \Phi(y_k)$ as above,
then $e_\alpha^* e_\alpha$ is also a c.a.i. for $\A$.  However
$e_\alpha^* e_\alpha  =  \sum_{k,l} \Phi(y_k)^* b_{kl} \Phi(y_l) $ 
where $b_{kl} = \Psi(x_k)^* \Psi(x_l) \in $ $\B$.  Since $B = [b_{kl}]$
is a positive matrix, it has a square root $R = [r_{ij}]$, say,
with entries $r_{ij} \in $ $\B$.  Thus $ e_\alpha^* e_\alpha  = 
\sum_k \Phi(y^\alpha_k)^* \Phi(y^\alpha_k) , $ where
$y^\alpha_k = \sum_j r_{kj} y_j$.   From this one can easily
deduce that the $\A-$valued innerproduct on $Y$ has dense
range, that is, $Y$ is a full right C$^*-$module over $\A$. 
Similar arguments show that $Y$ is a full left
C$^*-$module over $\B$, and that $X$ is also full on both 
sides.  Thus $X$ and $Y$ are strong Morita equivalence bimodules,
giving the strong Morita equivalence of $\A$ and $\B$.  

Observe that by the basic theory of
strong Morita equivalence (see e.g \cite{Ri2})
$_{\A}\K(X) \cong \B$.  Thus if $\{ f_{\beta} \}$ is a c.a.i.
 for $\B$, then $\{ G_{\beta} \}$ is a c.a.i.
 for $_{\A}\K(X)$, where $G_{\beta}(x) =
x f_{\beta} = G(r_{f_{\beta}})(x)$.   Observe too, by Lemma \ref{coron},
that
$F(W) \cong \{ T \in $ $_{\B}CB(B,F(W)) : T \; r_{f_{\beta}} \rightarrow T$
in norm $\}$ completely isometrically, where $\{ f_{\beta} \}$ is an
approximate identity for $\B$.  Applying the functor $G$ and Lemma
\ref{bf}, we see the last set is completely isometrically isomorphic to
$\{ S \in $ $_{\A}CB(X,GF(W)) : S G(r_{f_{\beta}}) \rightarrow S$
in norm $\}$, which is completely isometrically isomorphic to
$\{ S \in $ $_{\A}CB(X,W) :  S G_{\beta} \rightarrow S$
in norm $\}$, which in turn equals $_{\A}\K(X,W)$,
since $G_{\beta} \in $ $_{\A}\K(X)$.  Thus we have shown that
$F(W) \cong $ $_{\A}\K(X,W)$ completely isometrically, and it is an easy
algebra check that this is also as
left $\B-$modules.    Setting $W = $ $\A$ gives
 $Y \cong $ $_{\A}\K(X,A)$, so that $Y \cong \bar{X}$.
It is easily checked that this last relation is as
bimodules too.  
In Theorem 3.10 in \cite{Bhmo}, we showed that 
$\bar{X} \otimes_{h\A} W \cong $ $_{\A}\K(X,W)$
completely isometrically.  Thus $F(W) \cong Y \otimes_{h\A} W$
completely isometrically and as $\B-$modules, for all
$W \in $ $_{\A}OMOD$.  Its an easy algebra check 
now that $F \cong $ $_{\A}\K(X,-) \cong
Y \otimes_{h\A} -$ as functors.
By symmetry, we get the matching statement for $G$.  
The last statement of Theorem \ref{FT}, about the mapping of 
subcategories, follows because $\otimes_{h\A}$ coincides with 
the interior tensor product on the subcategories concerned.
\end{proof}

\vspace{7mm}

Thus the proof of our main result 
has boiled down to verifying the very concrete hypotheses of the
last theorem.
To that end, we first observe that the natural transformations
$GF(H) \cong H$ and
$FG(K) \cong K$  imply certain norm equalities.  Using, 
repeatedly, Lemmas \ref{coron}, \ref{bf} and the natural 
transformations, we see that
$$
\begin{array}{ccl}
H & \cong & \; _{\A}CB(\A,H) \cong \; _{\B}CB(Y,F(H))   
\; \cong \; _{\B}CB(Y,_{\B}CB(\B,F(H))) \\
& \cong & \; _{\B}CB(Y,_{\A}CB(X,GF(H))) 
\cong \; _{\B}CB(Y,_{\A}CB(X,H))
\end{array}
$$
completely isometrically.  Untangling these identifications shows that
$\zeta \in H$ corresponds to the following map $T_\zeta$ in the last space 
$_{\B}CB(Y,_{\A}CB(X,H))$ in 
the string above: namely $T_\zeta(y)(x) = (x,y)\zeta$.  
Thus 
$$
\begin{array}{ccl}
\Vert \zeta \Vert & = &  \Vert T_\zeta \Vert_{cb} \\
& = & \sup \{ \Vert [(x_{kl},y_{ij}) \zeta ] \Vert :
[x_{kl}] \in Ball(M_m(X)), [y_{ij}] \in Ball(M_n(Y)) , n,m \in \N \} \\
& = & \sup \{ \Vert [ \Psi(x_{kl}) \Phi(y_{ij}) \zeta ] \Vert :
[x_{kl}] \in Ball(M_m(X)), [y_{ij}] \in Ball(M_n(Y)) , n,m \in \N \} \\
& \leq & \sup \{ \Vert [ \Phi(y_{ij}) \zeta ] \Vert :
[y_{ij}] \in Ball(M_n(Y)) , n \in \N \} \\
& \leq & \Vert \zeta \Vert
\end{array}
$$
using equation (1), and the fact that $\Phi$ and $\Psi$ are complete 
contractions (Lemma \ref{phci}).
Thus $\Vert \zeta \Vert = \sup \{ \Vert [ \Phi(y_{ij}) \zeta ] \Vert :
[y_{ij}] \in Ball(M_n(Y)) , n \in \N \} $.  Squaring and using the usual 
formula for the matrix norms on $H^c$ we see that
$$
\begin{array}{ccl}
\langle \zeta \; \vert \; \zeta \rangle & = & \sup \{ \Vert [ 
\Phi(y_{ij}) \zeta ] \Vert^2 : [y_{ij}] \in Ball(M_n(Y)) ,
 n \in \N \} \\
& = & \sup \{ \Vert [ \sum_{k=1}^n \langle \Phi(y_{kj})\zeta \; 
\vert \; \Phi(y_{ki})\zeta 
\rangle ]\Vert : [y_{ij}] \in Ball(M_n(Y)) , n \in \N \} \\
& = & \sup \{ \Vert [ \langle (\sum_{k=1}^n 
\Phi(y_{ki})^*\Phi(y_{kj}))\zeta \;
\vert \; \zeta \rangle ]\Vert : [y_{ij}] \in Ball(M_n(Y)) , 
n \in \N \} \\ 
& = & \sup \{  \langle (\sum_{k=1}^n 
\Phi(\sum_{i=1}^n y_{ki}z_i)^*\Phi(\sum_{j=1}^n y_{kj}z_j))\zeta \;
\vert \; \zeta \rangle : [y_{ij}] \in Ball(M_n(Y)) , \sum_{i=1}^n
\vert z_i \vert^2 \leq 1 \} \\
\end{array}
$$
where the $z_i \in \C$.  Letting $y_k = \sum_{i=1}^n y_{ki}z_i$ we see 
that
$$\langle \zeta \; \vert \; \zeta \rangle = 
\sup \{  \langle (\sum_{k=1}^n \Phi(y_k)^* \Phi(y_k))\zeta \;
\vert \; \zeta \rangle  : [y_1, \cdots y_n]^t \in Ball(C_n(Y)) , 
n \in \N \} \; \; .  \eqno{(3)}   
$$
Replacing $\zeta$ by $\Psi(x) \eta$ for $x \in X, \eta \in K$ we have
$$
\langle \; \Psi(x)^* \Psi(x) \eta \; \vert \; \eta \rangle =
\sup \{  \langle (\sum_{k=1}^n \Psi(x)^*\Phi(y_k)^*\Phi(y_k)
\Psi(x))  \eta \; \vert \; 
\eta \rangle : [y_1, \cdots y_n]^t \in Ball(C_n(Y)) , 
n \in \N \} \; \; .  
$$
The expression $\sum_{k=1}^n \Psi(x)^*\Phi(y_k)^*\Phi(y_k)
\Psi(x)$ is, by equation (1) and the remark after it,
 an element $b \in \B$, with $0 \leq b \leq \Psi(x)^* \Psi(x)$ since 
$\Phi$ is completely contractive. 
Thus, for $x \in X, \eta \in K$ we have
$$
\langle \; \Psi(x)^* \Psi(x) \eta \; \vert \; \eta \rangle =
\sup \{  \langle \; b \eta \; \vert \; \eta \; \rangle : 
b \in \B \; , \; 0 \leq b \leq \Psi(x)^* \Psi(x) \} \eqno{(4)}
$$
A similar argument shows that  for $y \in Y, \zeta \in H$, we have
 $$
\langle \; \Phi(y)^* \Phi(y) \zeta \; \vert \; \zeta \rangle =
\sup \{  \langle \; a  \zeta \; \vert \; \zeta \rangle : 
a \in \A \; , \; 0 \leq a \leq \Phi(y)^*\Phi(y)  \} \eqno{(5)}
$$
It follows from (5), and the fact that every quasistate of $\A$ has a
unique w*-continuous extension to $e(\A)$ of
form $\langle \; \cdot \;  \zeta \; \vert \; \zeta \rangle$ for some
$\zeta \in Ball(H)$, that  $\Phi(y)^* \Phi(y)$ is a lowersemicontinuous
element in $e(\A) = \A''$ , for each $y \in Y$.  We refer the
reader to \cite{Ped} for details about lowersemicontinuity in 
the enveloping von Neumann algebra of a C$^*-$algebra.    
A similar, but slightly more complicated 
argument, shows that $\Psi(x)^* \Psi(x)$, as an element in   
$\B''$, corresponds to a lowersemicontinuous element in $e(\B)$
(which we recall, is W$^*-$isomorphic to $\B''$).  The complication
occurs since it seems we can say only that the 
quasistates of $\B$ have unique w*-continuous extensions to
$\B''$ of form $\sum_{k=1}^\infty \langle \; \cdot 
\; \eta_k \; \vert \; \eta_k \rangle$, where $\sum_{k=1}^\infty 
\Vert \eta_k \Vert^2 \leq 1 \; , \; \eta_k \in K$.   Nonetheless,
the calculation leading to equation (4) may be repeated, but
with $H$ and $K$ replaced by $H^\infty$ and $K^\infty$ 
(that is, the Hilbert space direct sum of a countably infinite number of 
copies of $H$ or $K$), to yield the desired conclusion.     

The crux of the proof now rests on a compactness argument
in $Q(\A)$, the (compact) set of quasistates of $\A$.  
For $y  = [y_1, \cdots y_n]^t 
\in Ball(C_n(Y))$, and $a_0 \in \A , \; 0 \leq a_0 \leq 1$, 
set $L_y = a_0 (\sum_{k=1}^n \Phi(y_k)^* \Phi(y_k)) a_0$,
which is a lowersemicontinuous element in $e(\A)$.  Moreover,
since $\Phi$ is completely contractive and since
$y \in Ball(C_n(Y))$, we see that $L_y \leq a_0^2$.   
Replacing $\zeta$ with $a_0 \zeta$ in  (3), and using
the fact that 
the quasistates of $\A$ are `vector quasi-states'
 of $e(\A)$,  we see that 
$$
\phi(a_0^2) = \sup \{ L_y(\phi) : y \in Ball(C_n(Y)) , n \in \N \} 
\eqno{(6)}  
$$
for all $\phi \in Q(\A)$.  Here $L_y(\phi)$ is the 
(scalar) value of $L_y$ (interpreted as an element of $\A^{**} = e(\A)$)
evaluated at $\phi \in $ $\A^*$.
 For $m \in \N$, and $y \in Ball(C_n(Y))$ set
$U_y^m = \{ \phi \in Q(\A) :  L_y(\phi) > \phi(a_0^2) 
(1 - \frac{1}{m}) - \frac{1}{m} \}$.
Since $L_y$ is lowersemicontinuous, $U_y^m$ is an open set in 
$Q(\A)$, and by (6) for each fixed  $m \in \N$, these sets 
$\{ U_y^m \}$ form an open cover of $Q(\A)$.   
Hence there is a finite subcover, corresponding to points
$y_1^m, \cdots , y_{k_m}^m$.

Keeping $m$ fixed, and $x \in BALL(X)$, we let 
$b^m_k = \Psi(x)^* L_{y_k^m} \Psi(x)$, which by equation
(1) and the remarks after it, is an element of $\B$.  Since
each $b^m_k$ is strictly
dominated (as a function on $Q(\B)$)
by the lowersemicontinuous function 
$\frac{1}{m} + \Psi(x)^* a_0^2 \Psi(x)$, it follows by a 
standard
lowersemicontinuity argument, effectively Dini's theorem
using \cite{Ped} Lemma 3.11.2, that there is an element $b^m \in \B$
satisfying $b^m \leq \frac{1}{m} + \Psi(x)^* a_0^2 \Psi(x)$, and also
$b^m \geq b^m_k - \frac{1}{m}$ for each $k$.    
It follows that for $\eta \in H, \Vert \eta \Vert = 1, $ and
$x \in BALL(X)$, that  
$$
\begin{array}{ccl}
\langle \; (\frac{1}{m} + \Psi(x)^* a_0^2 \Psi(x)) \eta \; \vert \; 
\eta \; \rangle
& \geq & \langle \; b^m  \eta \; \vert \; \eta \; \rangle \\
& \geq & \max_k  \langle \; b^m_k  \eta \; \vert \; \eta \; \rangle - \frac{1}{m} \\
& = & \max_k \;  \langle \; \Psi(x)^* L_{y_k^m} \Psi(x) 
 \eta \; \vert \; \eta \; \rangle - \frac{1}{m} \\
& = & \max_k \;  L_{y_k^m}(\phi_0)  - \frac{1}{m}  
\end{array}
$$
where $\phi_0(a) = \langle \; a \; \Psi(x) \eta \; \vert \; \Psi(x) \eta
\; \rangle$.  Since $x$ and $\eta$ have norm $\leq 1$, 
 $\phi_0$ is a quasistate.
Thus by the finite subcovering property we conclude that
$$
\begin{array}{ccl}
\langle \; ( \frac{1}{m} +  \Psi(x)^* a_0^2 \Psi(x)) 
\eta \; \vert \; \eta \; \rangle & \geq & \langle 
\; b^m  \eta \; \vert \; \eta \; \rangle \\
& \geq & \phi_0(a_0^2)(1 - \frac{1}{m}) - \frac{2}{m} \\
& = & \langle \; a_0^2 \; \Psi(x) \eta \; \vert \; \Psi(x) \eta
\; \rangle (1 - \frac{1}{m}) - \frac{2}{m}  
\end{array}
$$
Thus 
$$-\frac{1}{m}  \leq
 \Psi(x)^* a_0^2 \Psi(x) - b^m \leq \frac{1}{m} 
\Phi(x)^*a_0^2\Phi(x) + \frac{2}{m} \leq \frac{3}{m}$$ 
which shows that $b^m \rightarrow \Psi(x)^* a_0^2 \Psi(x)$ in norm.
Thus $\Psi(x)^* a_0^2 \Psi(x) \in $ $\B$.  Taking $a_0$ to be
element $e_\alpha$ in a c.a.i. for 
$\A$, shows that $\Psi(e_\alpha x)^*\Psi(e_\alpha x) \in \B$.
Thus $\Psi(x)^*\Psi(x) \in \B$.

A similar argument (which is slightly complicated by the fact that 
a quasistate of $\B$ is of the form $\sum_{k=1}^\infty \langle \cdot \;  
\eta_k , \eta_k \rangle$), 
shows that $\Phi(y)^* \Phi(y) \in \A$ for $y \in Y$, which by
Theorem \ref{slide} completes the proof.

{\bf Acknowledgements and Addenda:}
We particularly thank Gilles Pisier and Christian le Merdy for
answering some questions related to operator spaces, and
Vern Paulsen
for various conversations, and in particular for a suggestion
to move a certain trick to an earlier part of the proof
(which considerably shortened the proof).
Also, many of the ideas which are crucial
to this paper come from our collaboration
\cite{BMP,BMP2} with Vern Paulsen and Paul Muhly.  We thank Gerd Pedersen
for a discussion on lowersemicontinuity,  thank others for
encouragement to persist with this project, and thank
 the referee for his suggestions.

After finishing this paper in May 1997, we were
informed that P. Ara had also
obtained a characterization of strong Morita equivalence in terms
of isomorphism of module categories \cite{Ara1,Ara2}.  
However Ara works within a
quite different category, namely all modules in the sense of 
pure algebra, both left and right sided. These modules are not 
over the C$^*-$algebras, but over their Pedersen ideals.  Also,
the conditions on his functors (described in \cite{Ara2})
are also quite different.  
It might be interesting to try to  combine the various ideas
from our two characterizations.  
In any case, there is certainly no duplication of results or 
methods.

\setcounter{section}{3}

\end{document}